\documentclass[a4paper,11pt,reqno]{amsart}

\usepackage{amssymb}
\usepackage{amsmath}
\usepackage{amsthm}
\usepackage{amsfonts}
\usepackage{bbm}
\usepackage{color}

\usepackage{bookmark}
\usepackage{hyperref}
\hypersetup{pdfstartview={FitH}}

\newcommand{\R}{\mathbb{R}}

\newtheorem{theorem}{Theorem}

\newtheorem{conj}[theorem]{Conjecture}
{\normalsize }

\DeclareFontFamily{U}{mathx}{\hyphenchar\font45}
\DeclareFontShape{U}{mathx}{m}{n}{
<5> <6> <7> <8> <9> <10>
<10.95> <12> <14.4> <17.28> <20.74> <24.88>
mathx10
}{}
\DeclareSymbolFont{mathx}{U}{mathx}{m}{n}
\DeclareFontSubstitution{U}{mathx}{m}{n}
\DeclareMathAccent{\widecheck}{0}{mathx}{"71}

\numberwithin{equation}{section}

\begin{document}
\title{Singular Brascamp-Lieb: a survey}

\author[P. Durcik]{Polona Durcik}
\address{Polona Durcik, California Institute of Technology, 1200 E California Blvd, Pasadena CA 91125, USA}
\email{durcik@caltech.edu}

\author[C. Thiele]{Christoph Thiele}
\address{Christoph Thiele, Mathematisches Institut, Universit\"at Bonn, Endenicher Allee 60, 53115 Bonn, Germany}
 \email{thiele@math.uni-bonn.de}

\date{\today}

\begin{abstract}
 We present an overview of results on multi-linear singular integrals in the broader context of 
  Brascamp-Lieb inequalities.  This elaborates a lecture given at the inspiring conference on Geometric Aspects of Harmonic Analysis at
  Cortona 2018 in 
  honor of Fulvio Ricci. 
\end{abstract}

\maketitle

\section{Brascamp-Lieb forms and inequalities}
The recently active area of Brascamp-Lieb inequalities focuses
on invariant multi-linear forms in functions on Euclidean spaces.
By the Schwartz kernel theorem, the multi-linear forms acting on $n$-tuples of
Schwartz functions $F_j$ on $\R^{k_j}$ continuously in each 
argument are exactly the ones that can be written as
$$\Lambda(F_1\otimes F_2\otimes \dots \otimes F_n)$$
with a unique tempered distribution $\Lambda$ on $\R^{k_1+\dots +k_n}$.

Brascamp-Lieb forms arise when the distribution $\Lambda$
specializes to  integration over an affine subspace of $\R^{k_1+\dots +k_n}$
with respect to an invariant measure,

\begin{equation*}
\int_{\R^{k_1+\dots +k_n}}\Big(\prod_{j=1}^n F_j(x_j) \Big)\delta(\Pi (x-z))\, dx,
\end{equation*}
where $x$ denotes a vector with components $x_j$, 
 $\Pi$ is a linear map whose ker, translated by 
the vector $z$,  is the affine  space of integration, and  $\delta$ is the Dirac delta measure
on the range of the map $\Pi$. Here we have called the zero set of a linear map the  ker rather than the kernel
of the map so as to distinguish it from an integral kernel  such as for example in the
Schwartz kernel theorem.

A change of variables equates this form with 
\begin{equation*}
\int_{\R^{k_1+\dots +k_n}}\Big(\prod_{j=1}^n F_j(x_j+z_j) \Big)\delta(\Pi x)\, dx, 
\end{equation*}
which is a Brascamp-Lieb form with integration over a linear space, acting on
translates of the functions $F_j$. 
Using such a reduction, we shall assume throughout this survey that the space of integration
is linear, unless stated otherwise:  
\begin{equation}\label{blform}
\int_{\R^{k_1+\dots +k_n}}\Big(\prod_{j=1}^n F_j(x_j) \Big)\delta(\Pi x)\, dx.
\end{equation}

A further change of variables, replacing  $x$ by $x-z$ with a 
vector $z$ in the ker of $\Pi$, shows an invariance of the 
Brascamp-Lieb integral under translation of the functions by amounts $z_j$.
Similarly, one observes a homogeneity of the form under simultaneous
dilations of the functions.

Using the Fourier transform, one may write for a Brascamp-Lieb form
$$\widehat{\Lambda}(\widehat{F}_1\otimes \widehat{F}_2\otimes \dots \otimes
\widehat{F}_n), $$ 
where $\widehat{\Lambda}$ is integration over
the orthogonal complement of the subspace of integration of $\Lambda$. If $\Pi$ in \eqref{blform} is
an orthogonal projection, we may write for the Fourier transform integral
\begin{equation}\label{blformfourier}
\int_{\R^{k_1+\dots +k_n}}\Big(\prod_{j=1}^n \widehat{F}_j(\xi_j) \Big)\delta((1-\Pi) (\xi))\, d\xi.
\end{equation}
This allows to
identify further invariances of the form under simultaneous translations of the
Fourier transforms of the functions. A translation of the Fourier transform of a function is the same as a modulation of the function itself:
$$M_\xi F(x)=F(x)e^{2\pi i x\cdot \xi}.$$
Up to scalar multiples, the multi-linear forms of Brascamp-Lieb type
are determined by their translation and modulation symmetries.

One may write the integral over the subspace also as a parameterized integral.
 Assume the subspace has dimension $m$ and let
 $$I: \R^m\to \R^{k_1+\dots +k_n}$$
 be a parameterization. Denote by $I_j$ the 
 composition of $I$ with the projection onto the $j$-th coordinate space $\R^{k_j}$.
  We may then write for \eqref{blform}, up to scalar multiple,  
 \begin{equation}\label{embedbl}
 \int_{\R^m} \Big ( \prod_{j=1}^n F_j(I_j x)\Big )\, dx .\end{equation}
 Writing each $F_j$ as Fourier integral, we obtain for \eqref{embedbl}
 \begin{align*}
 & \int_{\R^{k_1+\dots +k_n}}\int_{\R^m} \Big ( \prod_{j=1}^n \widehat{F}_j(\xi_j)e^{2\pi i \xi_j \cdot (I_j x)}\Big )\, dx d\xi\\
 & =  \int_{\R^{k_1+\dots +k_n}}  \Big ( \prod_{j=1}^n \widehat{F}_j(\xi_j) \Big) \delta( \sum_{j=1}^n I_j^T \xi_j) d\xi ,
 \end{align*}
which is of the form  \eqref{blformfourier} with $1-\Pi=\sum_{j=1}^n I_j^T$.

It is natural to seek bounds for Brascamp-Lieb forms by products of norms of the functions,
 with a choice of norms respecting the symmetries of the form. Most common   are 
 Lebesgue  norms $\textup{L}^p$, which are invariant
 under translations and modulations and have a homogeneity under dilations. 
The corresponding  bounds  are called Brascamp-Lieb inequalities.
 With a choice of exponents  $p_j$, these inequalities are written as  
 \begin{equation}\label{brasliebineq}
  \Big | \int_{\R^m} \Big ( \prod_{j=1}^n F_j(I_j x)\Big )\, dx \Big | 
 \leq C \prod_{j=1}^n \|F_j\|_{p_j} 
  \end{equation}
 with a constant $C$ depending on the $I_j$ and $p_j$ but not on the Schwartz functions $F_j$.

Given a tuple of exponents, if $p_j<\infty$ for some $j$, then a Brascamp-Lieb inequality can only hold if the  map $I_j$ is surjective.
To see this, assume $I_j$ is not surjective. Let $y$ and $z$ parameterize respectively the range of $I_j$ and the orthogonal complement of this range in $\R^{k_j}$.
Then left-hand side of the Brascamp-Lieb inequality does not change under replacing $F_j$ by
$$\widetilde{F}_j(y,z):=F_j(y,\lambda z),$$ while the right-hand side scales with a power of $\lambda$ that
is non-trivial if $p_j<\infty$.

 If $p_j=\infty$, then the map $I_j$ need not be surjective. For example,  if $m=0$, then the 
projection $I_j$ is not surjective except in the the pathological case $k_j=0$. Nevertheless,
as the Brascamp-Lieb integral becomes evaluation 
at a point, the Brascamp-Lieb inequality holds with all exponents equal to $\infty$. 

Well known cases of a Brascamp-Lieb inequality are H\"older's inequality, where
all maps $I_j$ are the identity map, Young's convolution inequality, and the Loomis-Whitney inequality where $m=n$, 
$k_j=n-1$ and the one dimensional kers of the maps $I_j$ span the full space $\R^m$.

Much research has been devoted to  Brascamp-Lieb and related inequalities, we refer to \cite{BL}, \cite{BCCT}, 
\cite{BBBF}, \cite{BBCF} and the references therein. In particular, \cite{BCCT} proves a necessary and sufficient dimensional condition for a Brascamp-Lieb inequality to hold, namely that
  \begin{equation}\label{bcct}
  \dim(V)\le \sum_{j=1}^n \frac 1{p_j} \dim(I_j V)
  \end{equation}
  for every subspace $V$ of $\R^m$, with equality if $V=\R^m$. The easy direction of this equivalence is necessity of \eqref{bcct}.
It is seen by testing
the Brascamp-Lieb inequality on suitable characteristic functions $F_j$, generating them as limits of Schwartz functions.
The supports of these functions are such that the integrand on the left-hand side of \eqref{brasliebineq} is nonzero on a disc, more precisely on a one-neighborhood in $\R^m$ of a large ball in $V$ of radius $R$. The left-hand side of the Brascamp-Lieb inequality grows in $R$ with the order $R^{\dim(V)}$. The suitable choice of the function $F_j$  is the characteristic function of the projection of the disc to $\R^{k_j}$. Its $\textup{L}^{p_j}$ norms grow with the order $R^{\dim(I_j(V))/p_j}$. Letting  $R$ tend to infinity, we obtain the lower bound of \eqref{bcct}. The equality in case $V=\R^m$ is obtained by using in addition small balls in $\R^m$.

Since $\dim(I_j  \R^m)\le \dim(\R^m)$,  inequality   \eqref{bcct} for $V=\R^m$ in case $m>0$ implies that
\begin{equation}\label{sumofreciprocals} 1\le \sum_{j=1}^n \frac{1}{p_j}.\end{equation}
When  equality holds in \eqref{sumofreciprocals}, then each map $I_j$ is injective on $\R^m$ and we  obtain
$\dim(I_j  V)= \dim(V)$ for all subspaces $V$ of $\R^m$. In this case, the condition \eqref{bcct}
for $V=\R^m$ automatically implies the condition for all subspaces of $\R^m$. Assuming all $I_j$ are surjective
as well, which is a mild assumption given the previous discussion, all $I_j$ are bijective. Reparameterizing the range of
each $I_j$, we may assume that each $I_j$ is the identity map and  thereby identify H\"older's inequality.

While it may be tempting to study \eqref{brasliebineq} with some   $0< p_j<1$, such estimates are easily seen to fail.
This is also reflected by \eqref{bcct}. Assume for example a Brascamp-Lieb inequality with $p_1<1$ and denote the ker of $I_1$ by $W$.
Then we obtain a contradiction by applying \eqref{bcct} twice:
$$m=k_1+\dim(W)\le 
k_1+\sum_{j=2}^n \frac 1{p_j}\dim(I_j W)<
\frac {k_1}{p_1}+\sum_{j=2}^n \frac {k_j}{p_j} =m .$$
The endpoint case $p_j=1$ reduces to Brascamp-Lieb inequalities of  fewer functions.
We show this in case $j=1$. By a weak limiting process, the Brascamp-Lieb inequality extends to   finite Borel measures
 in place of the first Schwartz function. In particular, one may insert  translates of the Dirac delta measure.
Conversely, bounds for the Brascamp-Lieb integral with translates of the Dirac delta measure as the first input imply by superposition the Brascamp-Lieb inequality for arbitrary Schwartz functions as first input.
The Brascamp-Lieb inequality with a translate of the Dirac delta measure can be written as
\begin{align*}
\int_{\R^{k_1+\dots +k_n}}\Big(\prod_{j=2}^n F_j(x_j+y_j) \Big) \delta(x_1-y_1)\delta(\Pi x)\, dx,
\end{align*}
which can be further written as
\begin{align*}
 \int_{\R^{k_2+\dots +k_n}}\Big (\prod_{j=2}^n F_j(x_j+y_j) \Big ) \delta(\Pi(y_1,x_2,\dots, x_n))  \, dx_2,\dots dx_n. 
\end{align*}
Note that the  range of the the restriction of $\Pi$ to fixed $y_1$ is the same as the range of $\Pi$ as a consequence of the assumption that $I_1$ is surjective.
The last display is again a Brascamp-Lieb integral with an affine linear space of integration and one input function less. Thus
we have shown the desired reduction.
 
 This observation in reverse allows to interpret the Dirac delta measure in the general 
Brascamp-Lieb form \eqref{blform} as coming from an $\textup{L}^1$ function. Thus \eqref{brasliebineq} is
 equivalent to the inequality
$$
\Big| \int_{\R^{k_1+\dots +k_n}}\Big (\prod_{j=1}^n F_j(x_j) \Big )F_{n+1}(\Pi x)\, dx  \Big|\le
 C \Big ( \prod_{j=1}^n \|F_j\|_{p_j} \Big )\|F_{n+1}\|_1.
 $$
The integral on the left hand side is again a Brascamp-Lieb form  \eqref{blform}, if written as
$$\int_{\R^{k_1+\dots +k_n+k_{n+1}}} \Big (\prod_{j=1}^{n+1} F_j(x_j) \Big ) \delta(x_{n+1}-\Pi (x_1,\dots ,x_n))\, dx. $$
Here the subspace of integration is the graph of a function in the first $n$ variables.

\section{Singular Brascamp-Lieb inequalities}
 
Coming to the main subject of this survey, one may ask whether a variant of the
Brascamp-Lieb inequality continues to hold if one inserts singular integral kernels
instead of finite measures into one or several input slots with $p_j=1$. 
Singular integral kernels  in general fail to be finite measures, but in many situations one retains 
inequalities thanks to cancellation between positive and negative parts of the  kernel. 
Examples of singular integral kernels arise from integrating a mean zero Schwartz function
over the group of dilations
\begin{equation}\label{dilationintegral}
K(t)=\lim_{N\to \infty} \int_0^N \lambda^k \phi(\lambda t )\,  \frac{d\lambda}\lambda,\quad 
\widehat{K}(\tau)=\int_0^\infty \widehat{\phi}\Big (\frac  \tau  \lambda \Big ) \frac {d\lambda}{\lambda}.
\end{equation} 
Such kernels are  homogeneous under dilations and smooth outside the origin.
They are in general not locally integrable near the origin, yet they are tempered distributions
in the sense that the limit in $N$  has to be executed after the pairing with a Schwartz function. Tempered distributions with such limits are called principal value
distributions.
 More generally, one may consider tempered distributions $K$ on $\R^k$ whose
Fourier transform $\widehat{K}$, called the multiplier associated with $K$, is a bounded
measurable function satisfying the symbol estimates
\begin{equation}\label{czsymbol}|\partial^\alpha \widehat{K}(\tau)|\le C |\tau|^{-|\alpha|}
\end{equation}
for some constant $C$, all $\tau\neq 0$ and all multi-indices $\alpha$ up to suitably large order.  
This condition is satisfied for the above homogeneous kernels. For much of
our survey it is sufficient to consider these homogeneous kernels.
The Dirac delta measure is a singular integral kernel, it can be written in the form
\eqref{dilationintegral} with a Schwartz function  of integral zero, and its Fourier transform is
a  constant function. A simple way to ensure that a Schwartz function has integral zero is to make it odd.
Many of the interesting features of the theory can already be seen when restricting to odd kernels.

We  write singular Brascamp-Lieb inequalities as
\begin{equation}\label{singbraslieb}
\Big | \int_{\R^{m}} \Big( \prod_{j=1}^{h} F_j(I_j x)\Big) \Big(\prod_{j=h+1}^{n} K_j(\Pi_j x)\Big)
\, dx\Big |  \leq C \prod_{j=1}^h \|F_j\|_{p_j}
\end{equation}
with singular integral kernels $K_j$ on $\R^{k_j}$ and surjective  
maps 
$$I_j, \Pi_j: \R^m\to \R^{k_j}.$$ The constant $C$ is assumed to be independent
of the functions $F_j$, and is assumed to depend on the  kernels $K_j$ only through the constant
 in \eqref{czsymbol} and the bound on the order of derivatives  in \eqref{czsymbol}. 
For smooth homogeneous kernels, the constant $C$ is controlled by some Schwartz norm
of the Schwartz function $\phi$ in \eqref{dilationintegral}.

As we ask a given singular Brascamp-Lieb inequality to hold for all choices of singular integral kernels, it needs
to hold for the special choice of a Dirac delta measure. In particular, the bound \eqref{singbraslieb} needs to hold when all 
kernels are the Dirac delta measure. Note that
$$\prod_{j=h+1}^{n} \delta(\Pi_j x)=\delta (\Pi_{h+1} x,\dots ,\Pi_{n}x),$$
where the Dirac delta measure on the right-hand side lives in dimension 
$k_{h+1}+\dots +k_{n}$. In order for the integral in \eqref{singbraslieb}
to be well defined,  we need   the map
\begin{equation*}
x\mapsto (\Pi_{h+1} x,\dots ,\Pi_{n}x)
\end{equation*}
to be surjective. We assume this surjectivity and choose variables
$$t=(t_{h+1},\dots, t_{n})$$
on the range of this map. Changing coordinates and choosing $y$ as vector of
coordinates for the joint ker
\begin{equation}\label{jointkernel}
W=\bigcap_{j=h+1}^{n}\ker \Pi_j\ ,
\end{equation}
we may rewrite the integral in \eqref{singbraslieb} as
\begin{equation}\label{sblform}
\int_{\R^{m}} \Big( \prod_{j=1}^{h} F_j(I_j (y,t))\Big) \Big(\prod_{j=h+1}^{n} K_j(t_j)\Big)
\, dy dt.
\end{equation} 
Thanks to these conventions, it is particularly easy to
reduce a singular integral by setting one kernel $K_j$ equal to the Dirac delta measure.
One removes this kernel  from \eqref{sblform}, sets the coordinate $t_j$ equal to
zero, and removes the integration over the variable $t_j$.

The class of singular integral kernels is invariant under dilation symmetries but not under translation or modulation symmetries. 
The translation  symmetries of the Brascamp-Lieb integral discussed after \eqref{blform} leave the singular Brascamp-Lieb
form invariant only if the components $z_j$ in the notation after \eqref{blform} are zero for $j>h$, that is those $j$ belonging to  kernels. 
An analogous observation holds for the modulation symmetries.

The mean zero condition on the Schwartz function in \eqref{dilationintegral}
is an important theme in singular integral theory.
To see necessity of the cancellation, consider a kernel  $K_{n}$ of the form
\eqref{dilationintegral} generated by a non-negative Schwartz function that is not
constant equal to zero, and assume there is only one  kernel or 
reduce the complexity by replacing the other  kernels by Dirac delta measures. 
Consider \eqref{singbraslieb} with characteristic functions 
$F_j$ of standard unit balls in the respective dimensions similarly to the proof of necessity of \eqref{bcct}.
The right-hand side of \eqref{singbraslieb} is finite.
The  integrand  on the left-hand side is equal to $K_{n}(t_n)$ for $y$ in a small ball
about the origin and $t_{n}$ in a small fixed interval around the origin.  
Uniformly in this ball in $y$,  the integral in $t_{n}$ tends to $\infty$ with $N$,
because the degree of homogeneity of the singular integral kernel is critical for integration.
Thus the left-hand side of \eqref{singbraslieb} is unbounded.

Singular Brascamp-Lieb inequalities have  seen much development in recent years, but the level of
understanding is far from establishing a general criterion mirroring  the condition \eqref{bcct}.
We present some necessary and some sufficient conditions.

A necessary condition for  \eqref{singbraslieb} can be obtained by specifying all $K_j$ as Dirac delta measures, yielding a reduced  Brascamp-Lieb inequality of lower order with integration over the joint ker   defined in \eqref{jointkernel}.
We obtain that $I_j$ needs to map $W$ onto $\R^{k_j}$ if $p_j<\infty$, and \eqref{bcct} for the reduced inequality gives the necessary condition
\begin{equation}\label{newbcct}
  \dim(V)  \le \sum_{j=1}^{h} \frac 1{p_j} \dim(I_j V)
\end{equation}
for all $V\subseteq W$, with equality if $V=W$.

Due to the importance of cancellation of the singular integral kernel, we
may obtain further necessary conditions for  \eqref{singbraslieb}, namely that
 \begin{equation}\label{localcond}
\ker I_1+  \ker \Pi_n=\R^m ,
\end{equation}
and similarly for other   indices by permutation of the Schwartz functions and   kernels. 
To see necessity, assume this condition is violated.
By reduction we may assume $h=1$.
Then there is a non-zero linear functional  $\lambda$ on $\R^m$ which vanishes on
$\ker I_1$ and on $\ker \Pi_n$. 
This functional   factors as
\begin{equation*}
\lambda(x)=\rho_1(I_1 x)=\rho(\Pi_n x)
\end{equation*}
for some suitable maps $\rho_1$, $\rho$.
Let $K_n$ be the kernel defined by \eqref{dilationintegral} with the Schwartz function
$e^{-|t|^2}\rho(t)$ and define for any tuple of Schwartz functions 
 $$\widetilde{F}_1 =F_1 \times  ({\rm sgn}\circ \rho). $$
We obtain 
$$\int_{\R^m} F_1(I_1x)\Big( \prod_{j=2}^{n-1} F_j(I_j x)\Big) 
|K_n|(\Pi _nx)\, dx
=\int_{\R^m}{\widetilde{F}}_1(I_1 x) \Big( \prod_{j=2}^{n-1} F_j(I_j x)\Big)  K_n(\Pi_n x)\, dx.$$
Approximating $\widetilde{F}_1$ by Schwartz functions and applying a hypothetical singular
Brascamp-Lieb inequality for the right-hand side, we obtain the same inequality for the
left-hand side, contradicting the  impossibility of the inequality for the non-negative kernel $|K_n|$. 

If $p_1=\infty$, we obtain another necessary condition for a singular Brascamp-Lieb
inequality, which we adapt from \cite{vk:tp}, namely 
$$  \bigcap_{j=2}^{n-1} \ker I_j \subseteq \ker I_1 \cup  \ker \Pi_n \ .$$
For assume this is not the case. 
Pick a vector $u$ which is in the space on the left-hand side but not in the space 
on the right-hand side.
There is a linear functional $\lambda_1$ that factors as $\lambda_1(x)=\rho_1(I_1 x)$
and is positive on $u$. Let $F_1=1_{+}\circ \rho_1$ with $1_+$ the characteristic function of the positive
half line. Let $F_j$ for $2\le j\le h$ be the characteristic function of the unit ball. 

There is also a linear functional $\lambda$ that factors as $\lambda(x)=\rho(\Pi_n x)$ and is positive on $u$.
Let $K_n$ be the homogeneous kernel \eqref{dilationintegral} generated by $e^{-|t|^2}\rho(t)$.
We split the singular Brascamp-Lieb integral \eqref{singbraslieb} by first integrating along lines parallel to $u$:
$$\int _{\ker (\lambda)}\int _\R F_1(I_1(x+su))  \Big (\prod_{j=2}^{n-1} F_j(I_j (x+su)) \Big )K_n(\Pi_n(x+su))\, dsdx$$
 
 The  middle factor in the integrand, the product over $j$, is independent of $s$ and equl to $1$ for $x$ in a small neighborhood of the origin.
The first factor is bounded,
$$F_1(I_1(x+su)) =1_+(\lambda_1(x)+s(\lambda_1(u)),$$
and for some sufficiently large $a$ it vanishes for $s<-a$ and is constant $1$ for $s>a$.
The third factor is positive for $s>0$.
Hence the integral over $s<-a$ vanishes, is a bounded number for $-a<x'<a$, and is plus infinity for $s>a$
and $x$ in a small neighborhood of the origin. Hence the singular Brascamp-Lieb integral is unbounded.

We come to some sufficient conditions for singular Brascamp-Lieb inequalities to hold.
If  one of the exponents $p_j$ is equal to $1$, we may reduce a singular Brascamp-Lieb 
inequality to one of lower complexity by the use of Dirac delta measures as discussed in 
the non-singular case. Validity of the reduced inequalities becomes a sufficient criterion
for validity of the original inequality.

If
\begin{equation}\label{hycondition}
1\le p_j \le 2
\end{equation}
for all $1\le j\le h$,  then it is useful to pass to the integral on the
Fourier transform side. If $\Pi$ in \eqref{blform} is an orthogonal projection, the
Fourier transform integral reads as
\begin{equation}\label{ftside}
\int_{\R^{k_1+\dots +k_n}} \Big (\prod_{j=1}^{h} \widehat{F}_j(\xi_j)\Big ) 
\Big (\prod_{j=h+1}^n \widehat{K}_j(\xi_j) \Big )  \delta((1-\Pi)\xi)d\xi .
\end{equation}
This is estimated by  a non-singular Brascamp-Lieb
inequality in the Fourier transforms of the functions, using that the multipliers $\widehat{K}_j$ are  functions in $\textup{L}^\infty$. Aiming at the dual exponents
 ${p_j}'=p_j/(p_j-1)$, we need the    condition $\eqref{bcct}$:
$$\dim (V) \le \sum_{j=1}^{h}  \frac 1{{p_j}'} \dim{(V_j)} ,$$
where $V$ is a subspace of $\ker(1-\Pi)$, $V_j$ is its projection onto
the $j$-th coordinate space,  and equality holds for $V$ equal to $\ker(1-\Pi)$.
We thus estimate \eqref{ftside} with the Brascamp-Lieb inequality by
$$\le C \prod_{j=1}^{h} \|\widehat{F}_j\|_{{p_j}'} \le C \prod_{j=1}^{h} \|F_j\|_{p_j}.$$
In the second inequality we have used the Hausdorff Young inequality, which is applicable 
by the assumption \eqref{hycondition}.
An interesting variant of this theme is to estimate a singular Brascamp-Lieb integral
by a mixed product of $\textup{L}^p$ norms of the functions and $\textup{L}^p$ norms of the Fourier transforms of the functions.
An instance of this has been studied in \cite{kes2}.

\section{Inequalities with one singular kernel and H\"older scaling}

As seen in the previous  section, when all exponents $p_j$ are at most  $2$, then one has a good
sufficient criterion for a singular Brascamp-Lieb inequality.  At the other end of the spectrum, 
when the $p_j$ are large, one finds the special case of H\"older scaling
\begin{equation*}
 \sum_{j=1}^{h} \frac{1}{p_j} =1,
\end{equation*}
where in an average sense the $p_j$ are as large as they can be. This is a heavily
studied case and we shall assume it throughout the rest of the survey.

Recall that in the H\"older case the condition \eqref{newbcct}  needs only to be checked for $V=W$.
Each map $I_j$ restricted to $W$ needs to be injective. Neglecting some trivial
extensions for $p_j=\infty$, we may also assume that this map is surjective for each $j$.
As a consequence, all $k_j$, $1\le j\le n-1$ are equal and in particular $k_j=k_1$ and
 $$m=k_1+ k_n.$$  The singular  Brascamp-Lieb integral may then be written as
  \begin{align*}
\int_{\R^{k_1}}\int_{\R^{k_n}}   \Big ( \prod_{j=1}^{n-1} F_j(A_j y +B_j t) \Big ) K_n(t)  dtdy,
\end{align*}
with matrices $A_j$ and $B_j$. Each of the  matrices $A_j$ has to be regular.
Changing $F_j$ by precomposing with the matrix $A_j$, we
may assume that all $A_j$ are equal to the identity matrix,
  \begin{align}\label{sbl}
\int_{\R^{k_1}}\int_{\R^{k_n}}   \Big ( \prod_{j=1}^{n-1} F_j(y+B_j t) \Big ) K_n(t)  dtdy.
\end{align}
Interchanging the order of integration so that  $y$ becomes the inner variable and replacing it by $y-B_1t$, we may
in addition assume that
$$B_1=0.$$
Writing each $F_j$ as Fourier integral we obtain for \eqref{sbl}
\begin{align*}
&\int_{\R^{(n-1)k_1}}
\int_{\R^{k_1}}\int_{\R^{k_n}}   \Big ( \prod_{j=1}^{n-1} \widehat{F}_j(\eta_j) e^{2\pi i \eta_j \cdot (y+B_j t)} \Big ) K_n(t)  dtdy d\eta_1 \dots d\eta_{n-1}\\
&=\int_{\R^{(n-1)k_1}:\eta_1+\dots +\eta_{n-1}=0}
\int_{\R^{k_n}}   \Big ( \prod_{j=1}^{n-1} \widehat{F}_j(\eta_j) \Big ) \widehat{K}_n(- \sum_j B_j^T \eta_j )  dt d\gamma , 
\end{align*}
where $d\gamma$ is the Lebesgue measure on the subspace $\eta_1+\dots +\eta_{n-1}=0$ in $\R^{(n-1)k_1}$.

We  look at small values of $n$. For $n=2$,
the singular Brascamp-Lieb
integral in the discussed variables becomes
$$\int_{\R^{k_1}}\int_{\R^{k_2}} F_1(y)K_2(t)\, dt dy.$$
Taking formally the Fourier transform, one obtains
$$ \widehat{F}_1(0) \widehat{K}_2(0),$$
which is undetermined by \eqref{czsymbol} and does not lead to an interesting theory.

The case $n=3$ describes bilinear forms which dualize to linear operators.
In the above coordinates, the singular Brascamp-Lieb integral can be written as
$$\int_{\R^{k_1}}\int_{\R^{k_3}} F_1(y)F_2(y+Bt) K_3(t)\, dt dy.$$
If $B$ is not injective, we may integrate the ker of $B$ first. This
integrates the singular integral kernel towards a lower dimensional kernel, reducing the problem to a similar problem where $B$ is injective. If $B$ is not surjective, we may split the integration over $y$ into integration
over the range of $B$ and the complement of the range. The integral over the range
is a similar singular Brascamp-Lieb with smaller dimension, which can be estimated first.
Subsequently, one can estimate the complementary integral by H\"older's inequality.
Hence we may assume without loss of generality that $B$ is regular. By changing variables and replacing the kernel
$K_3$ by its composition with the inverse of $B$, we obtain the form
$$\int_{\R^{k_1}}\int_{\R^{k_3}} F_1(y)F_2(y+t) K_3(t)\, dt dy.$$
The dual linear operator is the classical convolution with a singular integral
kernel, which is well understood. As a consequence, we have the desired
singular Brascamp-Lieb inequality with H\"older scaling and
$1< p_1,p_2<\infty$.
The restriction $1<p_j$ can be understood as a condition of the type \eqref{localcond}
after a reduction by a Dirac delta function as in the discussion after \eqref{localcond}.

We turn to the genuinely multi-linear case $n\geq 4$. Fixing $n$ and $k_1$,
 singular Brascamp-Lieb inequalities become easier with growing $k_n$. 
 In case of odd kernels this can be made rigorous by the method of rotations, which
 we will discuss more thoroughly later.
   
The largest and thus easiest interesting case is $k_n=(n-2)k_1$. Beyond that, one would necessarily violate condition
\eqref{localcond} or be able to integrate out some of the $t$ variables of $K$ to reduce to a  kernel of smaller dimension.
The case $k_n=(n-2)k_1$ is the classical
theory of multi-linear operators of Coifman-Meyer type \cite{cm:czmo}.
Note that the map $(B_2 \otimes \ldots \otimes B_n)$ has to be surjective
or else one could again reduce the problem by integrating a trivial ker variable.
Changing coordinates to  parameterizing the range of this map  and  
adjusting the kernel $K_n$ suitably, we obtain
$$\int_{\R^{k_1}} \int_{\R^{(n-2)k_1}} F_1(y)\Big (\prod_{j=2}^{n-1} F_j(y+t_j)\Big ) K_n(t_2,\dots t_{n-1}) \, d(t_2,\dots, t_{n-1})dy.$$
With a further change of variables we may write more symmetrically
$$\int_{\R^{k_1}} \int_{\R^{(n-1)k_1}:\,t_1+\dots + t_{n-1}=0} \Big (\prod_{j=1}^{n-1} F_j(y+t_j)\Big ) \widetilde{K}_n(t_1,\dots t_{n-1}) \, d \gamma dy$$
with $d\gamma$ the invariant measure on the subspace of $\R^{(n-1)k_1}$ perpendicular
to the diagonal $(1,\dots, 1)$ and $\widetilde{K}_n$ suitably defined on this subspace.
As a result of the classical theory, one obtains singular Brascamp-Lieb inequalities with H\"older scaling  
as long as
$$1<p_j \le \infty$$
for all indices $1\le j\le n-1$. The restriction $1<p_j$ is again a consequence of the discussion after \eqref{localcond}. There is no restriction at $\infty$. An interesting theory allows to push
the inequalities of Coifman-Meyer type beyond infinity. 
Under certain conditions on the kernel, one obtains $BMO$ bounds, and one may consider restricted type estimates 
as discussed in \cite{mtt:mo}, dualizing bounds in earlier work \cite{ks:me}, \cite{gt:msi}.
Taking the Fourier transform, the Coifman-Meyer multi-linear form becomes 
$$\int_{\R^{(n-1)k_1}:\, \xi_1+\dots +\xi_{n-1}=0} \Big (\prod_{j=1}^{n-1} F_j(\xi_j)\Big ) \widehat{\widetilde{K}}_n(\xi_1,\dots \xi_{n-1}) \, d\gamma\, $$
where the Fourier transform of $\widetilde{K}_n$ is suitably taken in the space $\Gamma$.
The subspace of integration has dimension $(n-2)k_1$, which is equal to the
dimension $k$ of the multiplier. As a consequence, there are no translations of this subspace which leave the 
multiplier invariant. Hence the Coifman-Meyer case does not exhibit modulation symmetries.
It relies on  classical Calder\'on-Zygmund techniques that are translation and dilation invariant.

As one lowers $k$ from the maximal interesting $(n-2)k_1$, one may no longer
uniquely determine the embedding map   $I$ up to change of coordinates.
The discussion bifurcates depending on the geometry of $I$,
and the classification of cases leads to quite elaborate linear algebraic questions. 
One case in every dimension  is distinguished as the generic position of these projections. It can be obtained 
almost surely by picking $I$ randomly with respect to suitable Gaussian probability measures.
The study of this generic situation has begun in the work on the bilinear Hilbert transform
\cite{lt} and \cite{gn}. In the case $k_1=1$, the best  sufficient
dimensional condition in the generic situation is \cite{mtt:mo}. 
In the notation
$$\int_{\R}\int_{\R^{k_n}} \Big (\prod_{j=1}^{n-1} F_j(y+B_j t)\Big ) K_n(t) dt dx, $$
the generic case is when each tuple of the linear functionals $B_j$ spans the maximal possible space.
One obtains the singular Brascamp-Lieb inequality  with H\"older scaling for all
$$1<p_j\le \infty$$
provided one has the dimensional condition 
\begin{equation}\label{halfcond}
  n-3 <  2k_n
\end{equation}
for any $n\ge 3$.
Unlike the Coifman-Meyer case, the generic singular Brascamp-Lieb integral  for $k<(n-2)k_1$ exhibits modulation symmetries. 
The proof of the above result employs a modulation invariant counterpart of Calder\'on-Zygmund
techniques called time-frequency analysis. This technique originates in the works of \cite{carleson}, \cite{fefferman} and was first applied to singular Brascamp-Lieb forms in the work \cite{lt} on the bilinear Hilbert transform. An approach to time-frequency analysis through outer measures was described in \cite{dt}. The principal value limit in \eqref{dilationintegral} in the context of time-frequency analysis and in particular the bilinear Hilbert transform is studied in
\cite{l:max}, \cite{dop1}, \cite{dop2}.

While the  time-frequency analysis in \cite{mtt:mo}   breaks down if the condition
\eqref{halfcond} is violated, it remains an open problem whether \eqref{halfcond}
is necessary for singular Brascamp-Lieb inequalities to hold.
Even under condition \eqref{halfcond}, interesting open questions remain
concerning the extension of singular Brascamp-Lieb inequalities to restricted type inequalities
beyond the threshold at $p_j=\infty$. This  is discussed in \cite{mtt:mo}, see also \cite{dipt} for a discussion
near the boundary of the range of exponents with known bounds.

The extension of the above result of \cite{mtt:mo} to $k_1>1$ is addressed
in \cite{dpt},  proving singular Brascamp-Lieb inequalities on the form 
$$\int_{\R^{k_1}}\int_{\R^{k_n}} \Big (\prod_{j=1}^{n-1} F_j(y+B_j t)\Big ) K_n(t) dt dy$$
assuming $B_j:\R^{k_n}\to \R^{k_j}$ are in generic position and 
\begin{equation}\label{genhalfcond}
  k_1(n-3) <  2k_n.
\end{equation}
If $k_n$ is an integer multiple of $k_1$, this follows rather quickly from the methods of \cite{mtt:mo}.
For the fractional multiple case, \cite{dpt} uses some additional arguments from additive combinatorics.
The authors restrict attention to the range $2<p_j\le \infty$. It is not known whether the restriction $2<p_j$ is necessary.

A partial explanation for the break down of modulation invariant
time-frequency analysis beyond \eqref{halfcond}, \eqref{genhalfcond}
is the occurrence of more general symmetries. For example, consider the case of the trilinear
Hilbert transform
\begin{equation*}
\int_{\R}\int_{\R} \Big(\prod_{j=1}^4 F_j(y+B_j t)\Big ) \frac 1{t} dt dy
\end{equation*}
with generic, that is pairwise different, numbers $B_j$. This form
exhibits a symmetry under quadratic modulation
$$Q_{\alpha_j}F_j(x)=F_j(x)e^{i\alpha_j x^2}$$
where the four numbers $\alpha_j$ are all non-zero and satisfy
$$\sum_j\alpha_j (y+B_j t)^2=0.$$
It would be interesting to find extensions of time-frequency analysis
that are  invariant under more general symmetries and address boundedness
of the trilinear Hilbert transform.
This starts with a solid understanding of the type of symmetries,
we refer to related work on inverse theorems for Gowers norm
\cite{gt:ig3} involving generalized quadratic phase functions
possibly relevant for the trilinear Hilbert transform and the more
general symmetries in \cite{gtz:ign}.
A variant of time frequency analysis under polynomial symmetries was developed
in \cite{lie1}, \cite{lie2}, \cite{pz:mpm}.
Additional symmetries may not be the only obstruction to
go beyond \eqref{halfcond}, because it is not clear that all cases beyond
\eqref{halfcond} exhibit additional symmetries.

Shrinking $k_n$ further, the minimal non-trivial case is $k_n=1$. 
The  distance to $k_1$ is maximized  if $k_1=n-1=h$. If $k_1$ is greater than or equal to $h$, then the vectors $B_j$, $2\le j\le h$ span a space of dimension less than $k_1$ and one may reduce to a singular Brascamp-Lieb integral of lower order as discussed in the case $n=3$. By the same token, if $k_1=h$, then these
vectors have to be linearly independent and thus a basis of $\R^{k_1}$.
Since all bases are equivalent up to change of variables,  one can write the singular Brascamp-Lieb
integral without loss of generality in symmetric form as
\begin{equation}\label{simplex}
\int_{\R^{h}}\Big( \prod_{j=1}^h F_j(x_1,\ldots, x_{j-1},x_{j+1},\ldots, x_h)\Big )\frac{1}{x_1+\ldots + x_h} dx.
\end{equation}
This form is called the simplex Hilbert form.  
Maybe the biggest challenge in the area is to understand whether this
form satisfies any singular Brascamp-Lieb inequalities.
By symmetry and interpolation techniques, the easiest bound to prove should
be the one with all exponents equal. We formulate this as a conjecture.
\begin{conj}\label{simplexconj}
  There exists a constant $C$  such that for all tuples of Schwartz functions $(F_j)_{j=1}^h$
  the form \eqref{simplex} is bounded by
  \begin{align*}
         C  \prod_{j=1}^h \|F_j\|_{h}.
\end{align*}
\end{conj} 
By the method of rotations, bounds for the simplex Hilbert form imply
bounds for many singular Brascamp-Lieb integrals, including for the multi-linear Hilbert transform, another
major open problem. Moreover, bounds for the simplex Hilbert form imply   bounds for 
the Carleson and polynomial Carleson operator
\begin{align*}
\int_{\R}f(x-t)e^{i(N_1(x)t+N_2(x)t^2+\ldots + N_{d}t^{d})}\frac{dt}{t},
\end{align*}
 which was for general $d$ studied in \cite{lie1}, \cite{lie2} and \cite{pz:mpm}.
 Partial progress on the simplex Hilbert form in the case $h=3$ can
  be found in  \cite{ktz:tht}, which in particular establishes the above conjectured bound in a dyadic model
  when one of the functions takes a special form. Further results concerning 
 truncations of the simplex Hilbert  form and effective bounds in the parameter of truncation
 are discussed in \cite{pz:splx} based on the approach in \cite{tt}, and in \cite{dkt}.

 Having discussed generic choices of $B_j$ in the spectrum from large $k_n$ to small $k_n$, we turn attention to some 
 of the phenomena arising when we do not ask the $B_j$ to be in generic positions.
We begin with the simplest case which displays some of the phenomena,
$$\int_{\R }\int_{\R} \Big (\prod_{j=1}^3 F_j(y+B_j t)\Big ) K_4(t) dt dy .$$
The generic case has three different real numbers
$B_j$, this is the classical bilinear Hilbert transform. All generic cases have
the same proof of Brascamp-Lieb bounds using time-frequency analysis.
If two values of $B_j$ are equal, the form changes its nature. One identifies
the pointwise product of two functions, and replacing the product by a new function
we obtain a singular Brascamp-Lieb integral with $n=3$. Applying the classical theory
without time-frequency analysis and then applying H\"older's inequality to resolve
the product proves $\textup{L}^p$ bounds in this degenerate situation.
The case that all three values of $B_j$ are equal is even further degenerate but
of no interest, it leads to the pointwise product of three functions together
with the indeterminate integral in case $n=2$.
If two of the values of $B_j$ approach each other, the first proof of the
bilinear Hilbert transform produced a growing constant in the singular
Brascamp-Lieb inequality. It was natural to seek uniform bounds,
which was achieved in a series of papers \cite{ct:uni}, \cite{li:ub}, \cite{gl:ub}, \cite{ot}, \cite{uw}
in the full H\"older range of exponents with $1< p_j\le \infty$.
Some of these results were generalized to uniform bounds  on other families of
singular Brascamp-Lieb integrals in \cite{mtt:ue}.

A more complicated classification of cases occurs for the two dimensional bilinear Hilbert transform
\begin{equation*}
\int_{\R^2}\int_{\R^2} \Big (\prod_{j=1}^3 F_j(y+B_j t)\Big ) K_4(t) dt dy,
\end{equation*}
a situation first considered by \cite{dt:2dbht} and then thoroughly discussed
in the PhD thesis  \cite{mw}.
The thesis classifies the possiblilities for the   parameters $B_1$, $B_2$, $B_3$ into nine cases.
Most cases can be normalized such that $B_1=0$ and $B_2=I$, leaving only $B=B_3$ as indetermined matrix,
which may be assumed to be in Jordan canonical form.
A trivial pointwise product occurs if $B=0$ or $B=I$, this results in a reduction of the complexity
of the integral as in the one dimensional case.
The case that all eigenvalues of $B$ are different from $0$ and $1$ is the generic case covered by previous results.
The case that one eigenvalue of $B$ is equal to $0$ or $1$ and the other eigenvalue is different from
$0$ and $1$ is an interesting hybrid case discussed in \cite{dt:2dbht}, likewise the case of a
non-trivial Jordan block with eigenvalue $0$ or $1$.
The case when $B$ has both $0$ and $1$ as eigenvalue
is called the twisted paraproduct and is an instance of the forms in Theorem \ref{mainthm}
below with $m=2$, albeit with the fourth function set constant equal to $1$.  

Only in one of the nine cases it is not known whether the singular Brascamp-Lieb inequality holds
at a nontrivial set of exponents. This is the case where the first columns of all three
matrices $B_1,B_2,B_3$ vanish, while the second columns  respectively are $(0,0)$, $(0,1)$, $(1,0)$. 
This case is a simplex Hilbert form discussed in the above
conjecture.  All remaining cases reduce to easier objects and are of lesser interest.
An abundance of questions concerning uniform bounds arise between these various cases.  
While the method of rotations would prove uniform bounds for odd kernels 
from Conjecture \ref{simplexconj}, lacking a proof of  the latter it may be of interest to study these uniform questions.

We turn to a class of Brascamp-Lieb integrals where the modulation symmetry group is spanned by rich modulations symmetries.
 A rich modulation symmetry is  a modulation symmetry which generalizes to arbitrary phase functions.
 For example the H\"older form
 $$\int_{\R} F_1(x)F_2(x)\, dx$$
 is invariant not only under replacing $F_1$ and $F_2$ by
 $M_\xi F_1$ and $M_{-\xi} F_2$ respectively, but also under replacing them by
 $$F_1(x)e^{i\phi(x)},\quad F_2(x)e^{-i\phi(x)}$$
 for arbitrary real phase functions $\phi$.
 If we consider each input function as a function in $k_1$ arguments,
 then one way that rich modulations symmetries occur is when slots
of  different functions share the same argument.

 We consider an example where each of the $k=k_1$ slots carries two possible
variables, making it $2k$ integration variables, which we denote as
$$(x_1^0,\ldots ,x_{k}^0,x_1^1,\ldots ,x_{k}^1)=x.$$
Each possibe combination of the variable occurs in one of the functions.
This requires $2^{k}$ input functions parameterized by the cube
$Q$, the set of all
$$j:\{1,2,\dots, k\}  \to \{0,1\}.$$ 
Consequently, for $j\in Q$, we have   
\begin{align*}
I_j x = (x_1^{j(1)},x_2^{j(2)},\ldots ,x_k^{j(k)}).
\end{align*}
We further consider a singular integral kernel $K$ in $\R^{k}$ 
and an arbitrary surjective $\Pi: \R^{2k}\to \R^{k}$. The Brascamp-Lieb integral
in question then writes as  
 \begin{equation}\label{cubicbl}
\int_{\R^{m}} \Big( \prod_{j\in Q} F_j(I_j x)\Big ) K(\Pi x)\, dx\,  .
\end{equation}
\begin{theorem} [from \cite{dut}]
\label{mainthm}
Given $k\ge 1$, the form  \eqref{cubicbl} satisfies a singular Brascamp-Lieb
inequality with $p_j=2^{k}$ for all $j\in Q$ if and only if for all $j$
\begin{equation}
\label{cubebcct}
  k=  \dim(I_j (\ker \Pi)) .
\end{equation}
\end{theorem}

The condition \eqref{cubebcct} is the specialization of \eqref{newbcct} in this situation.

 While rich symmetries are very large symmetry groups and restrict
 techniques to those that are invariant under these symmetries, at least they
 have a very generic structure and one does not need to delve into the theory
 of polynomial or other structured symmetries.  The main technique 
 in the context of rich symmetries  was pioneered  in
 the context of the so-called twisted paraproduct in  \cite{vk:tp} and is sometimes called twisted technology. 
Brascamp-Lieb integrals involving rich symmetries were also studied in \cite{vk:bell},  \cite{bernicot:fw}, \cite{kt:t1},  \cite{pd:L4}, \cite{pd:Lp} and also in  \cite{dkst}, \cite{kas} with applications to quantitative convergence of ergodic averages, and  in \cite{dkr}, \cite{patterns}  with applications to some problems in Euclidean Ramsey theory. An application to stochastic integrals was studied in \cite{ks}. Further higher dimensional generalizations are discussed in   \cite{ms}.

It would be desirable to study some natural extensions of  Theorem \ref{mainthm}. One obvious generalization would be a more general range
of exponents than the symmetric exponent point.
Somewhat related to that is the question what happens if the 
corners of the cube are not fully occupied, that is the number of functions
is strictly less than $2^{k_1}$. In case one has $\textup{L}^\infty$ bounds, it is trivial
to omit the corresponding function by estimating the constant function
in $\textup{L}^\infty$, but it is not clear that all inequalities with constant
functions arise from more general $\textup{L}^\infty$ bounds.

One further extension is to allow more than two variables in one slot, that is 
for $k_1\ge 1$ and $l\geq 2$ we may consider  $\R^{m}$ with coordinates
$$x=((x_1^0,\ldots, x_{k}^0), (x_1^1,\ldots, x_{k}^1),\ldots, (x_{1}^{l-1},\ldots, x_{k}^{l-1}))\in \R^{k l}$$
 Then for all $j:\{1,2,\ldots, k\}\rightarrow \{0,\ldots,l-1\}$ we may define  
$$I_{j} x = (x_1^{j(1)},x_2^{j(2)},\ldots, x_{k}^{j(k)})$$
One may then ask the analoguous result as Theorem \ref{mainthm}.  

Note that also the simplex Hilbert forms of Conjecture \ref{simplexconj}
have many rich modulation symmetries. Indeed, the group of modulation symmetries 
of the simplex Hilbert form is spanned by  rich symmetries. The space of integration
in Fourier space has dimension $n(n-2)+1$. Since the singular integral kernel is one dimensional, this gives
 an $n(n-2)$ dimensional group of modulation symmetries of the simplex Hilbert form.
 However, for each of the $n$ variables one can find  $n-2$ pairs of functions so that independent rich symmetries akin to the above shown
 apply between this pair of functions.
The forms in Theorem \ref{mainthm} and the suggested generalization above 
have the structure that each variable has a fixed slot number in which it
may occur. Note that this is not the case in the simplex Hilbert form.
 For example, the variable $x_2$ typically appears in the second slot,
unless in the function $F_1$, where the variable $x_1$ is omitted and
the variable $x_2$ appears in the first slot. This mismatch is the main obstacle to
apply twisted technology to the simplex Hilbert form.

\section{Method of rotations and more general kernels}

 The method of rotation allows to write a  singular Brascamp-Lieb form with one singular integral kernel as a
 superposition of a family of forms with lower dimensional  kernels. The family of forms  is generated
by rotations or more general linear transformations of the space of integration.

Turning to details, a  singular Brascamp-Lieb form with a homogeneous smooth kernel can be written as
\begin{equation}\label{bl1}
\int _0^\infty  \int_{\R^{k_1+\dots +k_{n} }} \Big (\prod_{j=1}^{n-1} F_j(x_j) \Big ) t^{k_{n}}\psi(tx_{n}) \delta(\Pi x)\, dx  \frac{dt}t
\end{equation}
with a smooth and compactly supported  function $\psi$ with integral zero.  Assume 
there is a vector $v$ such that the inner product $v\cdot x_{n}$ is bounded away from zero on the support of $\psi$.
The following display is a superposition by a weight function $\phi$ of a family of forms generated by rank one perturbations of  $\Pi$ using a further fixed vector 
$w$ and a varying scalar parameter $a$:  
 \begin{equation*}
 \int_\R \int _0^\infty    \int_{\R^{k_1+\dots +k_{n} }}\Big (\prod_{j=1}^{n-1} F_j(x_j) \Big ) 
 t^{k_{n}} \psi (tx_{n})  \phi(a) \delta(\Pi x+w a (v\cdot x_{n})) \, dx \frac{dt}t da.
 \end{equation*}
We assume $\phi$ is smooth and compactly supported.
Rescaling the variable $a$ and combining it with the vector $x_n$ to a vector of dimension $k_{n}+1$, 
we recognize a new singular Brascamp-Lieb form
\begin{equation}\label{bl2}
\int _0^\infty     \int_{\R^{k_1+\dots +(k_{n} +1)}}\Big (\prod_{j=1}^{n-1} F_j(x_j)\Big ) 
t^{k_{n}+1}\widetilde{\psi} (tx_{n}, ta)  
 \delta(\Pi x+  w a ) \, dx da \frac{dt}t
 \end{equation}
with the compactly supported smooth function
\begin{equation}\label{decomp}
\widetilde{\psi}(x_{n},a):=\psi (x_{n})  
\frac 1{|v\cdot x_{n}|}\phi \Big (\frac{a}{v\cdot x_{n}}\Big ).
\end{equation}
One verifies that $\widetilde{\psi}$ has integral zero by first integrating in  $a$ and then in $x_n$.
 If we can prove bounds for the singular Brascamp-Lieb forms  \eqref{bl1}
uniformly for all maps $\Pi$ in the perturbed family, then  by superposition
we obtain a bound with the same exponents for \eqref{bl2}.

Conversely, given a Brascamp-Lieb integral  as in \eqref{bl2}, one may seek
to write it as superposition of Brascamp-Lieb forms with lower dimensional kernels.
A general procedure exists, when the function $\widetilde{\psi}$ is odd.
In addition, we assume $\widetilde{\psi}$ is compactly supported away from the origin.
After  a  decomposition  by a finite smooth partition of unity,
and a suitable rotation of the coordinate system  for each piece, we can assume that there is a vector 
$v$ of dimension $k_n$ such that $\widetilde{\psi}$ is supported in the union of two small neighborhoods 
respectively of $(v,0)$ and $(-v,0)$

With suitable compactly supported functions $\varphi$ and $\rho$
we may write 
$$\widetilde{\psi}(x_{n},a)
=\frac 1{|v\cdot x_{n}|}\varphi\Big (x_{n},\frac a{v\cdot x_{n}} \Big ) 
=\frac 1{|v\cdot x_{n}|}\varphi\Big (x_{n},\frac a{v\cdot x_{n}} \Big ) \rho\Big (\frac a{v\cdot x_{n}} \Big )
 $$
and note that $\varphi$ is odd in the  first variable for fixed second variable.
Taking a Fourier integral  of $\varphi$ in the second variable and denoting that by $\widehat{\varphi}$, we obtain
$$\widetilde{\psi}(x_{n},a)
=\int_\R \widehat{\varphi}(x_{n},\xi)
\frac 1{|v\cdot x_{n}|}
e^{2\pi i \xi \frac a{v\cdot x_{n}}} \rho\Big (\frac a{v\cdot x_{n}} \Big )\, d\xi .
 $$
 For fixed $\xi$, the integrand is a function of the form \eqref{decomp} with an odd function 
$\psi$. If we can prove bounds for the family of Brascamp-Lieb integrals
of lower dimensional kernels uniformly for fixed Schwartz norm of $\psi$ of some order,
then we may integrate these bounds in $\xi$ as the Schwarz norm of $\widehat{\varphi}$ 
in the first variable is rapidly decreasing as a function in the second variable.

One can iterate rank one perturbations to obtain the more general superposition
$$
 \int_{\R^l} \int _0^\infty    \int_{\R^{k_1+\dots +k_{n} }} \Big (\prod_{j=1}^{n-1} F_j(x_j) \Big ) 
 t^{k_{n}} \psi (tx_{n})  \phi(a) \delta \Big (\Pi x+ \sum_{i=1}^l w_i a_i (v_i\cdot x_{n}) \Big ) \, dx \frac{dt}t da .
$$
If the function $\phi$ in the above calculation is replaced by a finite Borel measure,
in particular a Dirac delta measure, estimates for the form \eqref{bl2} are equivalent
to estimates for the form \eqref{bl1} with lower dimensional kernel uniformly
over the perturbation parameters in the support of $\phi$. Choosing $\phi$ with
any intermediate regularity between smooth function and Borel measure, one can view
the difficulty of estimates for the superposed operator  as intermediate between
the two endpoint cases. Estimates for such forms with rough singular integral kernel 
can be of their own interest, if estimates for the lower dimensional kernels are not known
or maybe known to be false in general.

An early example of this principle   is provided by the 
Calder\'on commutator \cite{calderon}, which  later appeared in the investigation  of  the Cauchy integral along Lipschitz curves, see  \cite{CMM82} and the references therein. 
The commutator can be viewed as a rough superposition of bilinear Hilbert transforms.
Calder\'on proposed the study of the bilinear Hilbert transform and uniform bounds for it  as a stepping stone towards  the  commutator. However, the bilinear Hilbert transform remained an open problem for many years after bounds for the  Calder\'on commutator were obtained using different techniques.  A recent account and  approach to the Calder\'on commutator  and higher order commutators  was given in \cite{Mus14a} and in \cite{Mus14b}. These higher order commutators can be seen as a suitable superposition of multi-linear Hilbert transforms which by themselves are not known to be bounded.

If $\Pi$ as in \eqref{blform} is perturbed by a rank one map, then the embedding map $I$ as in \eqref{embedbl} 
can also be identified  as perturbed by a rank one matrix.  To be more precise, we assume that the perturbation
is $\Pi+\Pi(u)\otimes v$ where $v$ is a vector in $\ker\Pi$ and $u$ is orthogonal to $\ker\Pi$.
This representation can be found if the perturbation is small and the dimension of the ker of the perturbed map is equal to that of the
original map, namely $m$, but the kers are different.  As we have a rank one perturbation,  the two kers
intersect in a space of dimension $m-1$, and we may choose a unit vector $v$ in $\ker\Pi$ perpendicular to this subspace.
Using that the perturbation is small, we may chose $u$ perpendicular to $\ker\Pi$ so that $v-u$ is in the ker
of the perturbation. Then $\Pi+\Pi(u)\otimes v$ has the same ker as the perturbation and we may assume it is the perturbation.
The perturbation of the embedding map $I$ can then be written as $I- u\otimes I^T v$. To verify this, one checks separately that the vectors
that embed under $I$ into the intersection of the kers of $\Pi$ have the same image under the perturbed map,
and that the vector that maps to $v$ under $I$ maps to $v-u$ under the perturbation.

If the perturbations are such that only one component $u_j$ of $u$ and only the component $(I^Tv)_n$ of $I^Tv$ is non-zero,  we may view the averaging
of the form as an averaging of the function $F_j$. If we iterate several perturbations like that, then the  averaged function  takes the form
$$ F(I_j x, x_n)=\int_{\R^{l}} \phi(a)  F_j \Big (I_j x- \Big (\sum_l a_l u_l  (v_l \cdot I x_n)\big )_j \Big )   \, da .$$
If there are enough averages so that the rank one matrices add to a regular matrix, and  if $F$ is in $\textup{L}^\infty$,  then the averaged function $F(y,z)$  becomes a $y$ dependent 
symbol in the variable $z$ in the sense
$$|\partial_y^\alpha  \partial_z^\beta   F(y,z)|  \le C |z |^{-|\alpha|-|\beta |}$$
 for all multi-indices up to some degree depending on the regularity of the averaging function $\phi$.
Multiplying  this symbol with the  singular integral kernel gives a "space dependent" singular integral
form which is nowadays seen  within in the theory of $T(1)$ theorems originating in \cite{dj}. Therefore, bounds
for the averaged operator can be viewed as a Brascamp-Lieb version of a $T(1)$ theorem.

In this spirit, a multi-linear $T(1)$ theorem with a  variant of   the bilinear Hilbert transform with space dependent singular integral kernel was proven in \cite{bdnttv}
and applied in \cite{pal} in a singular variant of a higher Calder\'on commutator. $T(1)$ theorems with rich modulation symmetries 
were proven in \cite{kt:t1}, \cite{ms} in dyadic models, it would be interesting to extend these results to the continuous setting
and extend to further averaged singular Brascamp-Lieb forms.

The paper \cite{DR18}  discusses  averages of the simplex Hilbert forms which yield singular Brascamp-Lieb 
forms with rich modulation symmetries. The averaged forms are such that they can be treated by twisted technology. More precisely,
\cite{DR18} proves bounds in cases $n=4$ and $n=5$ on
\begin{align*}
\int_{(0,1)^{n-3}}\int_{\R^{n-2}}\int_{\R}  \Big ( \prod_{j=1}^{n-3} F_j(y+\alpha_j B_j t) \Big ) F_{n-2}(y+B_{n-1}t)F_{n-1}(y)
K_n(t)  dtdy d\alpha
\end{align*}
for linearly independent vectors $B_j$.

\section{Inequalities with two singular kernels and H\"older scaling}

Singular Brascamp-Lieb integrals in the case of several singular integral kernels fall into the scope
of multi-parameter theory. We display some of the features of multi-parameter theory
using the example of two kernels. We continue to assume H\"older scaling.

 Considerations  analoguous to those leading to \eqref{sbl} from \eqref{sblform} turn the singular Brascamp-Lieb 
 integral with two kernels into the form
 \begin{equation}\label{twosbl}
\int_{\R^{k_1}}\int_{\R^{k_{n-1}}} \int_{\R^{k_{n}}}
 \Big(\prod_{j=1}^{n-2} F_j(y +B_js +C_jt) \Big ) K_{n-1}(s)K_n(t) dtdsdy  .
\end{equation}
 Applying the Fourier transform as after \eqref{sbl} we obtain the alternative expression
 \begin{equation}\label{twofourier}
 \int_{\Gamma}
\Big ( \prod_{j=1}^{n-2} \widehat{F}_j(\xi_j) \Big ) 
\widehat{K}_{n-1}(- \sum_{j=1}^{n-2} B_j^T \xi_j )  
\widehat{K}_n(- \sum_{j=1}^{n-2} C_j^T \xi_j ) 
\, d\gamma , 
\end{equation}
where $\Gamma$ is the subspace of $\R^{(n-2)k_1}$ determined by 
$\xi_1+\dots +\xi_{n-2}=0$ and $d\gamma$ is the Lebesgue measure on this subspace.

 Simplifying degenerations may occur. The arguments of the two multipliers in \eqref{twofourier} can be identical,
 that is each $C_j$ is equal to $B_j$. As the product of two multipliers is again a multiplier with analoguous symbol bounds, this reduces
 to a singular Brascamp-Lieb with one kernel. 
 Another simplifying degeneration of \eqref{twosbl} may be separation. If for every $j$ one of the matrices $B_j$ or $C_j$
 is zero, then we may write the integral  in $s$ and $t$ as a product of two integrals, one in
 $s$ and one in $t$. Then we may apply H\"older's inequality in the variable $x$ on this product.
Resolving  the resulting $\textup{L}^p$ norms by pairing with a dual function, we obtain  
two singular Brascamp-Lieb integrals with one kernel each.
 Separation in \eqref{twosbl}  may occur after replacing the variable $y$ by $y+Bs+Ct$ for suitable
 matrices $B$ and $C$.

 A family of cases 
 occurs with counterexamples to a singular Brascamp-Lieb inequality
that show a phenomenon not possible for one kernel. Assume we have a family of quadratic forms $Q_j$ on $\R^{k_1}$ such that
$$\sum_{j=1}^{n-2} Q_j(y+B_js+C_jt)=s_1 t_1 $$
where $s_1$ and $t_1$ are the first components of $s$ and $t$, there being no
loss in generality choosing these particular components. 
For $n$ large enough compared to $k_1, k_{n-1},k_n$, such quadratic forms will exist in the case of generic 
matrices $B_j$ and $C_j$. Choose functions of the form
$$F_j(x)=\phi(x) e^{-2\pi iQ_j(x)}$$
where $\phi$ is a non-negative smooth approximation of the characteristic function of a very large ball about the origin.
Choose the kernel
$$K_n(t)=\lim_{N\to \infty}\int_0^N \lambda^{k_n} \psi(\lambda t_1) \phi(\lambda (t_2,\dots,t_{k_n})) \frac {d\lambda}{\lambda}$$ 
with odd $\psi$ which is non-negative  on the positive half axis and with non-negative $\phi$, and similarly
for $K_{n-1}$ with  odd $\widetilde{\psi}$ such that $\widehat{\widetilde{\psi}}=\psi$.
Zooming into the critical integrals in $s_1$ and $t_1$ in the expression \eqref{twosbl}, we see
$$\int_\R \int_\R  e^{-2\pi i s_1t_1}\widetilde{\psi}(\mu s_1) \psi (\lambda t_1) \, ds_1  dt_1
=\mu^{-1}  {\psi}(\mu^{-1}t_1) \psi(\lambda t_1).$$
 The right-hand side is an even function in $t_1$ and non-negative on the positive half axis, hence it is 
 non-negative, and it is not identically zero as one can see considering
  $\mu^{-1} $ near $\lambda$.
 The effect is that the cancellation of the kernel $K_n$ is destroyed, 
 resulting in unboundedness as $N$ tends to $\infty$.  More details of this calculation
 can be found in  \cite{mptt:bi} for the two examples
 \begin{equation*}
   \int_{\R^4}F_1(x_1,x_2) F_2(x_1-t, x_2-s) F_3(x_1+t, x_2+s)
   \frac{ds}{s}\frac{dt}{t} dx
 \end{equation*}
 and 
 \begin{align}\label{secondquadratic}
\int_{\R^3}F_1(x) F_2(x+t) F_3(x+s)F_4(x+t+s) \frac{ds}{s}\frac{dt}{t} dx.
 \end{align}

  Multi-parameter theory is named after the various scaling parameters occurring in a product of
  singular integral kernels. We call the product of the multipliers in
  \eqref{twofourier} the joint multiplier and write it with scaling parameters $\mu$ and $\lambda$ as 
  \begin{equation*}
  m(\sigma ,\tau)=\widehat{K}_{n-1}(\sigma )\widehat{K}_{n}(\tau )=
  \lim_{N,M\to \infty} \int_0^N \int_0^M  
  \widehat{\phi}_{n-1}\Big (\frac{\sigma }{\mu}\Big ) \widehat{\phi}_n\Big (\frac{\tau }{\lambda}\Big )\  \, \frac{d\mu}{\mu} \frac{d\lambda}{\lambda}.
  \end{equation*}
  A typical step in multi-parameter theory is the cone decomposition, which is a sorting of an integral in several scaling parameters by the size
  of the scaling parameters as follows:
   \begin{equation*}
    m_1(\sigma ,\tau)+m_2(\sigma ,\tau)=\lim_{N\to \infty} \int_{0< \mu<\lambda<N}  \dots   \, \frac{d\mu}{\mu}\frac{d\lambda}{\lambda}+
   \lim_{M\to \infty} \int_{0<\lambda<\mu<M}  \dots  \, \frac{d\mu}{\mu}\frac{d\lambda}{\lambda}.
   \end{equation*}
  Note that the  joint multiplier $m$  in \eqref{twofourier} satisfies the multi-parameter symbol estimate
  \begin{equation}\label{firstpartials}
    |\partial_\sigma^\alpha \partial_\tau^{\beta} m(\sigma ,\tau)|\le C|\sigma |^{-|\alpha|}
    |\tau |^{-|\beta|},
  \end{equation}
  where $\partial_\sigma$ and $\partial_\tau$ are any partial derivatives in the $\sigma$ and $\tau$ variables respectively.
  The cone multipliers $m_1$ and $m_2$ satisfy
   \begin{align}\label{firstcone}
 & |\partial_\sigma^\alpha  \partial_\tau^{\beta} m_1(\sigma ,\tau)|\le C|\sigma|^{-|\alpha|-|\beta|}\  ,\\
\label{secondcone}
 & | \partial_\sigma ^\alpha \partial_\tau^{\beta}  m_2(\sigma, \tau)|\le C|\tau|^{-|\alpha|-|\beta|}\ .
   \end{align}
In some instances, bounds for the variants of \eqref{twofourier} with
the joint multiplier replaced by the cone multipliers can be established,
based on the symbol estimates \eqref{firstcone}, \eqref{secondcone}.
Note that these symbol estimates, say \eqref{firstcone}, are generalizations of the single kernel case $K_{n-1}=\delta$
in that the multiplier \eqref{firstcone} is "frequency dependent" in the variable
$\tau$, a dual concept to the "space dependent" kernels discussed in the previous section.
Typically, estimates for the cones hold for generic choices of the matrices $B_j$ and $C_j$
provided the methods of \cite{mtt:mo} or \cite{dpt} for "frequency dependent" multipliers 
apply, which is under the suitably adapted conditions \eqref{halfcond} and \eqref{genhalfcond}.
An example for a singular Brascamp-Lieb form where this cone decomposition applies and uses generalized bounds for
"frequency dependent" variants of the bilinear Hilbert transform is given by
\begin{align*}
\int_{\R^3} F_1(y)F_2(y+s+t)F_3(y+B_3s+C_3t)K_4(s)K_5(t) {ds} {dt} dy
\end{align*}
with generic parameters $B_3$ and $C_3$.

Somewhat opposite of the case of generic matrices $B$, $C$, one finds in the literature the case when each of these matrices is
either zero or elementary, meaning it  has precisely one non-zero entry, and this entry is equal to one.
The flag paraproducts in \cite{Mus:flag1}, \cite{Mus:flag2} are 
essentially this case for $k_1=1$. Estimates are shown for the case
$$\int_{\R^5}F_1(y)F_2(y-t_1) F_3(y-t_2-s_1) F_4(y-s_2)K_5(t_1,t_2)K_6(s_1,s_2)ds_1ds_2dt_1dt_2dy,$$
which is motivated by questions in fluid dynamics, and a rather general positive conjecture is formulated  in \cite{Mus:flag1}, \cite{Mus:flag2}.
While one also does a cone decomposition in this case, it is important that the multiplier retains a product structure underneath the cone decomposition, and one does not simply rely on symbol estimates \eqref{firstpartials}.
Necessity of the product structure is demonstrated in \cite{gk}. While a form
$$\int_{\R^3}F_1(y)F_2(y+t) F_3(y+s) )K_4(t)K_5(s)ds dt dy$$
is bounded by the method of separation, and the joint multiplier satisfies
$$|\partial_\sigma^{\alpha} \partial_\tau^{\beta}
(\widehat{K_4}(\sigma)\widehat{K_5}(\tau))|\le C|\sigma|^{-\alpha}|\tau|^{-\beta},$$
the form obtained by replacing the joint multiplier by a general multiplier $m$
satisfying
$$|\partial_\sigma^{\alpha} \partial_\tau^{\beta}
m(\sigma, \tau)|\le C|\sigma|^{-\alpha}|\tau|^{-\beta}$$
need not satisfy any bounds in $\textup{L}^p$ spaces.

We may consider the case of $B_j$ and $C_j$ being zero or elementary for $k_1>1$ as well. A particular instance is discussed in \cite{mptt:bi} under the name of bi-parameter paraproduct:
\begin{align*}
  \int_{\R^6} F_1(y_1,y_2)F_2(y_1+s_1, y_2+t_1)F_3(y_1+s_2, y_2+t_2)
  K_4(s_1,s_2)K_5(t_1,t_2) dsdtdy.
\end{align*}
A generalization with more   kernels is discussed in \cite{mptt:multi}.
These examples are not affected
by the obstruction described in \cite{gk}, and one may prove bounds for
multipliers satisfying \eqref{firstpartials}.
However, already a simple modification of the above such as interchanging
$s_2$ and $t_2$ in the argument of $F_3$ is not addressed by the discussion in \cite{mptt:bi}.

A hybrid between the generic case and the flag paraproduct case is called
the biest and studied in  \cite{mtt:biest1}, \cite{mtt:biest2},
 \begin{align*}
\int_{\R^3}F_1(x) F_2(x+t) F_3(x+s)F_4(x-t-s) \frac{ds}{s}\frac{dt}{t} dx.
 \end{align*}
 It arises in the theory of iterated Fourier integrals, which occur
 in multi-linear expansions of certain ordinary differential equations.
 Singular Brascamp-Lieb inequalities for this form are known and require time frequency analysis because the bilinear Hilbert transform is embedded into this object.
Compare with the similar form \eqref{secondquadratic}.
For a study of objects related to the biest see \cite{mtt:splx}, \cite{kes1}, \cite{kes2}, \cite{kes3}, \cite{jj}, \cite {dmt}.

A more recent development is the theory of vector valued inequalities in the
context of singular Brascamp-Lieb inequalities. 
The helicoidal method was introduced in \cite{BM1} to study forms similar to the biest through mixed norm spaces and vector-valued inequalities. A survey of the helicoidal method can be found in \cite{BM2}.

\section{Acknowledgements}

This survey was initiated during a delightful stay at the conference \emph{Geometric Aspects of Harmonic Analysis} in honor of Fulvio Ricci 2018 in Cortona, Italy.
The second author acknowledges support by the Deutsche Forschungsgemeinschaft
through the Hausdorff Center for Mathematics, DFG-EXC 2047, and
the Collaborative Research Center 1060.


\begin{thebibliography}{03}
\bibitem{BM1}
C. Benea, C. Muscalu, \emph{Multiple vector-valued inequalities via the helicoidal method}. Anal. PDE, 9(8):1931–1988, 2016.

\bibitem{BM2}
C. Benea, C. Muscalu, \emph{The helicoidal method}. Preprint (2018), arXiv:1801.10071.



\bibitem{BCCT} J. Bennett, A. Carbery, F. M. Christ, and T. Tao,
\emph{The {B}rascamp-{L}ieb inequalities: finiteness, structure and extremals.} Geom. Funct. Anal. {\bf 17} (2008), no.~5, 1343--1415. 
 
 \bibitem{BBBF} J. Bennett, N. Bez, S. Buschenhenke, T. C. Flock,
 \emph{The nonlinear Brascamp-Lieb inequality for simple data}.
 Preprint (2018), arxiv:1801.05214.
 
 \bibitem{BBCF}	J. Bennett, N. Bez, M. G. Cowling, T. C. Flock, 
 \emph{Behaviour of the Brascamp-Lieb constant}.
 Bull. Lond. Math. Soc.
{\bf 49} (2017),
no.~3, 512--518.




 
\bibitem{bdnttv}	
 A. Benyi, C. Demeter, A. Nahmod, R. Torres, C. Thiele F. Villarroya, \emph{Modulation invariant bilinear $T(1)$ theorem}. J. Anal. Math. {\bf 109} (2009), 279--352.
 
 
  \bibitem{bernicot:fw}
F.~Bernicot, \emph{Fiber-wise Calder\'on-Zygmund decoposition and application to a bi-dimensional paraproduct}. Illinois J. Math.  {\bf 56} (2012), no. 2, 415-422. 


 \bibitem{BL}
 	H. J. Brascamp, E. Lieb,
 	 	\emph{Best constants in Young's inequality, its converse, and its
 		generalization to more than three functions}.
 	Advances in Math.
 	{\bf 20} (1976),
 	no.~2, 151--173.
 	
\bibitem{calderon}
A.-P. Calder\'on, 
\emph{Commutators of singular integral operators.}
Proc. Nat. Acad. Sci. U.S.A. {\bf 53} (1965), 1092--1099.  	

\bibitem{carleson}
L. Carleson,
\emph{On convergence and growth of partial sums of Fourier series.}
Acta Math. {\bf 116} (1966), 135--157. 

  	\bibitem{CMM82} R. R. Coifman, A. McIntosh, Y. Meyer.
	\newblock L'int\'{e}grale de Cauchy d\'efinit un op\'erateur born\'e sur $L^2$ pour les courbes lipschitziennes.
	\newblock {\em Ann. of Math. (2)} {\bf 116(2)} (1982), 361--387.
	
	
      \bibitem{cm:czmo}
        R. Coifman, Y. Meyer, \emph{Calder\'on-Zygmund and mulilinear operators. Translated from the 1990 and 1991 French originals by David Salinger.} Cambridge Studies in Advanced Mathematics {\bf 48} Cambridge University Press, Cambridge, 1997.


\bibitem{dj}
 David, G., Journ\'e, J.-L., \emph{A boundedness criterion for generalized Calder\'on-Zygmund operators.} Ann. of Math. (2) {\bf 120} (1984), no. 2, 371--397.


        \bibitem{dpt} C. Demeter, M. Pramanik, C. Thiele, \emph{Multilinear singular operators with fractional rank.}
          Pacific J. Math. {\bf 246} (2010), no. 2, 293--324.
          
	\bibitem{dt:2dbht}
	C. Demeter, C. Thiele, \emph{On the two-dimensional bilinear Hilbert transform}. Amer. J. Math. {\bf 132} (2010), no. 1, 201--256.

      \bibitem{dipt}
        F. Di Plinio, C. Thiele, \emph{Endpoint bounds for the bilinear Hilbert transform.}
         Trans. Amer. Math. Soc. {\bf 368} (2016), no. 6, 3931--3972. 
        
       \bibitem{dmt}
         Y. Do, C. Muscalu, C. Thiele,
         \emph{Variational estimates for the bilinear iterated Fourier integral.} J. Funct. Anal. {\bf 272} (2017), no. 5, 2176--2233. 

       \bibitem{dop1}
         Y. Do, R. Oberlin, E. Palsson,
         \emph{Variational bounds for a dyadic model of the bilinear Hilbert transform.}
         Illinois J. Math. {\bf 57} (2013), no. 1, 105--119. 

       \bibitem{dop2}
         Y. Do, R. Oberlin, E. Palsson,
\emph{Variation-norm and fluctuation estimates for ergodic bilinear averages.} Indiana Univ. Math. J. {\bf 66} (2017), no. 1, 55--99. 

         
       \bibitem{dt}
	Y. Do, C. Thiele,
\emph{Lp theory for outer measures and two themes of Lennart Carleson united.} 
Bull. Amer. Math. Soc. (N.S.) {\bf 52} (2015), no. 2, 249--296. 




       \bibitem{pd:L4}
	P. Durcik, \emph{An $L^4$ estimate for a singular entangled quadrilinear form}. Math. Res. Lett. {\bf 22} (2015), no. 5, 1317-1332.
	
	\bibitem{pd:Lp}
	P. Durcik, \emph{$L^p$ estimates for a singular entangled quadrilinear form}. Trans. Amer. Math. Soc. {\bf 369} (2017), no. 10, 6935-6951.
 	
 	\bibitem{patterns}
P. Durcik, V. Kova\v{c}, 
 \emph{Boxes, extended boxes, and sets of positive upper density in the Euclidean space}. Preprint (2018), arXiv:1809.08692. 	
 
 
  	\bibitem{dkr}
P. Durcik, V. Kova\v{c}, L. Rimani{\'c}, 
 \emph{On side-lengths of corners in positive density subsets of the Euclidean space}.  Int. Math. Res. Not. 2018, no. 22, 6844-6869.
 
 	\bibitem{dkst} P. Durcik, V. Kova\v{c}, K. \v{S}kreb, C. Thiele, 
 	\emph{Norm-variation of ergodic averages with respect to two commuting transformations.} Ergodic Theory Dynam. Systems {
 	\bf 39} (2019), no. 3, 658--688.
 	
\bibitem{dkt} P. Durcik, V. Kova\v{c},  C. Thiele, 
 	\emph{Power-type cancellation for the simplex Hilbert  transform.}
To appear in J. Anal. Math.

\bibitem{DR18}P. Durcik, J. Roos, \emph{Averages of simplex Hilbert transforms}. Preprint (2018), arXiv:1812.11701.
 
\bibitem{dut} P. Durcik,  C. Thiele, 
 	\emph{Singular Brascamp-Lieb inequalities.}
Preprint (2018), arXiv:1809.08688.


\bibitem{fefferman}
 C.  Fefferman
\emph{Pointwise convergence of Fourier series.}
Ann. of Math. (2) {\bf 98} (1973), 551--571. 

\bibitem{gn} J. Gilbert, A. Nahmod, \emph{Bilinear operators with non-smooth symbols. I.}
  J. Fourier Anal. Appl. {\bf 7} (2001), no. 5, 435--467.

\bibitem{gk} L, Grafakos, N. Kalton, 
  \emph{The Marcinkiewicz multiplier condition for bilinear operators.}
   Studia Math. {\bf 146} (2001), no. 2, 115--156. 

  


\bibitem{gl:ub}L. Grafakos, X. Li \emph{Uniform bounds for the bilinear Hilbert transforms. I.}
Ann. of Math. (2) {\bf 159} (2004), no. 3, 889--933. 



\bibitem{gt:msi}L. Grafakos, R. Torres \emph{On multilinear singular integrals of Calderón-Zygmund type.} 
 Publ. Mat. 2002, Vol. Extra, 57--91. 

\bibitem{gt:ig3}
 B. Green, T. Tao,\emph{
   An inverse theorem for the Gowers $U^3(G)$ norm.}
 Proc. Edinb. Math. Soc. (2) {\bf 51} (2008), no. 1, 73--153. 

\bibitem{gtz:ign}
   B. Green, T. Tao, T. Ziegler, \emph{
  An inverse theorem for the Gowers $U^{s+1}[N]$-norm.}
Ann. of Math. (2) {\bf 176} (2012), no. 2, 1231--1372. 


\bibitem{jj}
J. Jung, \emph{Iterated trilinear Fourier integrals with arbitrary symbols}. Preprint (2013), arXiv:1311.1574.
 
\bibitem{ks:me} C. Kenig, E. Stein,\emph{Multilinear estimates and fractional integration.}
Proceedings of the 6th International Conference on Harmonic Analysis and Partial Differential Equations (El Escorial, 2000).
 Math. Res. Lett. {\bf 6} (1999), no. 1, 1--15.
 
 
\bibitem{kes1}
R. Kessler, \emph{Generic Multilinear Multipliers Associated to Degenerate Simplexes}. Collectanea Mathematica. doi:10.1007/s13348-018-0224-z. 

\bibitem{kes2} 
R. Kessler, \emph{Mixed Estimates for Degenerate Multilinear Operators Associated to Simplexes}. J. Math. Anal. Appl. {\bf 424} (2015) 344–360. 

\bibitem{kes3} 
R. Kessler, \emph{$L^p$ Estimates for Semi-Degenerate Simplex Multipliers}. Preprint (2016), arXiv:1609.05964.
 

	\bibitem{vk:bell}
	V. Kova\v{c}, \emph{Bellman function technique for multilinear estimates and an application to generalized paraproducts}. Indiana Univ. Math. J. {\bf 60} (2011), no. 3, 813--846.
	
	\bibitem{vk:tp}
	V. Kova\v{c}, \emph{Boundedness of the twisted paraproduct}. Rev. Mat. Iberoam. {\bf 28} (2012), no. 4, 1143--1164.
	
	
	\bibitem{ks}
	V. Kova\v{c}, K. A. \v{S}kreb, \emph{One modification of the martingale transform and its applications to paraproducts and stochastic integrals}. J. Math. Anal. Appl. {\bf 426} (2015), no. 2, 1143-1163.
	
\bibitem{kt:t1}
 Kova\v{c}, V., Thiele, C., \emph{A $T(1)$ theorem for entangled multilinear dyadic Calder\'on-Zygmund operators.} Illinois J. Math. {\bf 57} (2013), no. 3, 775--799. 	
	
	
	\bibitem{ktz:tht}
	V. Kova\v{c}, C. Thiele, P. Zorin-Kranich, \emph{Dyadic triangular Hilbert transform of two general and one not too general function}. Forum of Mathematics, Sigma {\bf 3} (2015), e25.
 	
\bibitem{l:max}
M. Lacey,
\emph{The bilinear maximal functions map into Lp for $2/3<p \le 1$}.
Ann. of Math. (2) {\bf 151} (2000), no. 1, 35--57. 

        
 \bibitem{lt} M. Lacey, C. Thiele,
 	\emph{$L^p$ estimates for the bilinear Hilbert transform}.
 	Proc. Nat. Acad. Sci. U.S.A.,
 	{\bf 94}
 	(1997),
 	no.~1,
 	33--35.

\bibitem{li:ub} X. Li \emph{ Uniform bounds for the bilinear Hilbert transforms. II. }
Rev. Mat. Iberoam. {\bf 22} (2006), no. 3, 1069--1126. 

\bibitem{lie1}
V. Lie, \emph{The (weak-$L^2$) boundedness of the quadratic Carleson operator}. Geom. Funct. Anal., {\bf 19.2} (2009), pp. 457–-497.

\bibitem{lie2}
V. Lie, \emph{The polynomial Carleson operator}. Preprint (2011), arXiv:1105.4504.


	\bibitem{Mus14a} C. Muscalu. 
	\newblock Calder\'on commutators and the Cauchy integral on Lipschitz curves revisited: I. First commutator and generalizations.
	\newblock {\em Rev. Mat. Iberoam.}, {\bf 30} (2014), 727--750.
	
	\bibitem{Mus14b} C. Muscalu. 
	\newblock Calder\'on commutators and the Cauchy integral on Lipschitz curves revisited II. The Cauchy integral and its generalizations.
	\newblock {\em Rev. Mat. Iberoam.}, {\bf 30} (2014), 1089--1122.
	
\bibitem{Mus:flag1}
C. Muscalu, \emph{Flag paraproducts}. Contemp. Math {\bf 505}, 131-151.

	
\bibitem{Mus:flag2}
C. Muscalu, \emph{Paraproducts with flag singularities I: A case study}. Revista Mat. Iberoamericana {\bf 23}, 705-742.


\bibitem{mptt:bi}
C. Muscalu, J. Pipher, T. Tao, C. Thiele,  \emph{Bi-parameter paraproducts}. Acta Math. {\bf 193} (2004), no. 2, 269--296.

\bibitem{mptt:multi}
C. Muscalu, J. Pipher, T. Tao, C. Thiele,  \emph{Multi-parameter paraproducts}. Rev. Mat. Iberoam. {\bf 22} (2006), no. 3, 963--976.


\bibitem{mtt:biest1}
 C. Muscalu, T. Tao, C. Thiele, \emph{$L^p$ estimates for the biest I. The Walsh case}.
Math. Ann. {\bf 329} (2004), no. 3, 401--426. 

\bibitem{mtt:biest2}
 C. Muscalu, T. Tao, C. Thiele, \emph{$L^p$ estimates for the biest II. The Fourier case}.
Math. Ann. {\bf 329} (2004), no. 3, 427--461. 

\bibitem{mtt:splx}
 C. Muscalu, T. Tao, C. Thiele, \emph{Multilinear operators associated to simplexes of arbitrary length}. Advances in analysis: the legacy of Elias M. Stein, 346–401, Princeton Math. Ser., 50, Princeton Univ. Press, Princeton, NJ, 2014.

\bibitem{mtt:mo} C. Muscalu, T. Tao, C. Thiele, \emph{Multi-linear operators given by singular multipliers}.
 J. Amer. Math. Soc. {\bf 15} (2002), no. 2, 469--496.
         
\bibitem{mtt:ue} C. Muscalu, T. Tao, C. Thiele, \emph{Uniform estimates on multi-linear operators with modulation symmetry. Dedicated to the memory of Tom Wolff.}
J. Anal. Math. {\bf 88} (2002), 255--309. 

    
	\bibitem{pal} E. Palsson.
\emph{$L^p$ estimates for a singular integral operator motivated by Calder\'{o}n's second commutator}.
 J. Funct. Anal., {\bf 262} (2012), 1645--1678.

\bibitem{ot} R. Oberlin, C. Thiele,
  \emph{New uniform bounds for a Walsh model of the bilinear Hilbert transform.}
  Indiana Univ. Math. J. {\bf 60} (2011), no. 5, 1693--1712. 



\bibitem{ms}
	M. Stip\v{c}i{\'c}, \emph{$T(1)$ theorem for dyadic singular integral forms associated with hypergraphs}. Preprint (2019), arXiv:1902.10462.

\bibitem{kas}
	K. A. \v{Skreb}, \emph{Norm-variation of cubic ergodic averages}. Preprint (2019), arXiv:1903.04370.     

 \bibitem{tt} T. Tao,
 	\emph{Cancellation for the multilinear Hilbert transform}.
 	Collect. Math. 
 	{\bf 67}
 	(2016),
 	no.~2,
 	191--206.
 

      \bibitem{ct:uni} C. Thiele, 
        \emph{A uniform estimate.}
       Ann. of Math. (2) {\bf 156} (2002), no. 2, 519--563. 


\bibitem{uw} G. Uraltsev, M. Warchalski, \emph{Uniform bounds for the bilinear Hilbert transform in local $L^1$.}
       Chapter in the PhD Thesis of M. Warchalski, RFWU Bonn, 2018.
       
      \bibitem{mw}
  M. Warchalski,
  \emph{Uniform estimates in one-and two-dimensional time-frequency analysis}.
  PhD Thesis, RFWU Bonn, 2018.
   
	\bibitem{pz:splx}
	P. Zorin-Kranich, \emph{Cancellation for the simplex Hilbert transform}. Math. Res. Lett. {\bf 24.2} (2017), pp. 581–592.
	
		\bibitem{pz:mpm}
	P. Zorin-Kranich, \emph{Maximal polynomial modulations of singular integrals}. Preprint (2017), arXiv:1711.03524.
	
\end{thebibliography}
\end{document}